\documentclass[11pt]{article}
\usepackage{authblk}
\usepackage{hyperref}
\usepackage{amsmath,amssymb, amscd}
\usepackage{amsfonts}
\usepackage{multicol}
\usepackage[latin1]{inputenc}
\usepackage{graphicx}
\usepackage{color}
\usepackage{cleveref}
\usepackage{fancybox}
\usepackage[]{fancyhdr}
\usepackage{enumerate}
\usepackage{makeidx}
\usepackage[all]{xy}
\usepackage{array}
\usepackage{calc}
\usepackage{indentfirst}
\usepackage{url}
\usepackage[shortlabels]{enumitem}
\usepackage{xcolor}
\usepackage{tikz}
\usetikzlibrary{snakes,arrows,shapes}
\newtheorem{example}{Example}
\newtheorem{theorem}{Theorem}
\newtheorem{lemma}[theorem]{Lemma}
\newtheorem{corollary}[theorem]{Corollary}
\newtheorem{remark}[theorem]{Remark}

\newtheorem{proposition}[theorem]{Proposition}
\newcommand{\prova}{\noindent{\bf Proof:\ }}
\begin{document}
\title{The spectrum of an $I$-graph}

\author[1]{Allana S. S. de Oliveira}
\author[2]{Cybele T. M. Vinagre}

\affil[1]{\small
Instituto de Ci\^{e}ncias Exatas, Universidade Federal Fluminense, Volta Redonda, Brasil. e-mail: allanasthel@id.uff.br}
 \affil[2]{\small
Instituto de Matem\'atica e Estat\'{\i}stica, Universidade Federal Fluminense, Niter\'{o}i, Brasil. e-mail: cybl@vm.uff.br}
 \maketitle
%
%\begin{frontmatter}
%
%\title{The spectrum of an $I$-graph}
%%\tnotetext[mytitlenote]{Fully documented templates are available in the elsarticle package on \href{http://www.ctan.org/tex-archive/macros/latex/contrib/elsarticle}{CTAN}.}
%
%%% Group authors per affiliation:
%%\author{Elsevier\fnref{myfootnote}}
%%\address{Radarweg 29, Amsterdam}
%%\fntext[myfootnote]{Since 1880.}
%
%%% or include affiliations in footnotes:
%\author[allanaaddress]{Allana S.S. de Oliveira}
%%\ead[url]{www.elsevier.com}
%
%\author[cyaddress]{Cybele T.M. Vinagre\corref{mycorrespondingauthor}}
%\cortext[mycorrespondingauthor]{Corresponding author}
%\ead{cybl@vm.uff.br}
%
%\address[allanaaddress]{Instituto de Ci\^{e}ncias Exatas, Universidade Federal Fluminense, Volta Redonda, Brasil}
%\address[cyaddress]{Instituto de Matem\'atica e Estat\'{\i}stica, Universidade Federal Fluminense, Niter\'oi, Brasil}

\begin{abstract}
We completely determine
the spectrum of an $I$-graph, that is, the eigenvalues of its adjacency matrix. We apply our result to prove known characterizations of connectedness and  bipartiteness in $I$-graphs by using an spectral approach. With our result, we also  determine the nullity of a certain subfamily of $I$-graphs.
%For integers $n\geq 3$ and $1\leq
%j,k<\frac{n}{2}$, the  $I$-graph $I(n, j, k)$ is the graph with vertex set
%$ \{a_0, a_1, \ldots, a_{n-1},  b_0,$ $ b_1, . . . , b_{n-1}\}$ \
%and edge set \
%$ \{\{a_i, a_{i+j}\},$ $ \{a_i, b_i\}, \{b_i, b_{i+k}\} : i = 0, \ldots , n-1\}$,  with subscripts reduced modulo $n$. We completely determine
%the spectrum of an $I$-graph, that is, the eigenvalues of its adjacency matrix. We apply our result to obtain two known properties
%of $I$-graphs by using spectral techniques and also, to determine the nullity of a certain subfamily.
\end{abstract}
\noindent \textbf{AMS Subject Classification}: 05C50, 15B05.\\
\noindent \textbf{Keywords}: $I$-graph, generalized Petersen graph, nullity.
%\begin{keyword}
%   Adjacency matrix\sep $I$-graph \sep generalized Petersen graph \sep nullity.
%\MSC[2010] 05C50\sep  15B05
%\end{keyword}

%\end{frontmatter}

%\linenumbers

\section{Introduction}

The class of $I$-graphs was  introduced in the Foster Census
\cite{origemigrafo}
 as  a natural
generalization of the so called \cite{origemnome} generalized Petersen graphs and has attracted the attention of many graph theorists.
Considerable study of $I$-graphs under algebraical,  combinatorial and geometric approaches can be found in       \cite{igrafo}, \cite{igrafohamiltoniano}, \cite{Horvatisom}, \cite{Petenum} and   references therein.
In our work we   investigate the $I$-graphs under an spectral approach.

In our work, we completely determine the eigenvalues of the adjacency matrix of an $I$-graph by using known properties of circulant and circulant block matrices. Furthermore, we apply our result and spectral graph techniques to
give new proofs of known  necessary and sufficient conditions for
  bipartiteness and connectedness of arbitrary $I$-graphs. Also,  we establish
the nullity, that is, the dimension of the eigenspace associated to the null eigenvalue, for a certain subfamily of $I$-graphs.

\section{Preliminaries}

The \emph{adjacency matrix} $\textbf{A}(G)=[a_{ij}]$ of an arbitrary simple graph $G$ with vertices $x_1, x_2, \ldots, x_n$, is the $n \times n$ matrix  where $a_{ij}=1$, if  $\{x_i, x_j\}$ is an edge of $G$, and $a_{ij}=0$ otherwise.
Since a simple graph has no loops or undirected edges, its adjacency matrix is  symmetric and has zero diagonal elements.  The
\emph{characteristic polynomial} $p_{G}(\lambda )$
of
$G$
is that of
$\textbf{A}(G)$, that is,
$p_{G}(\lambda )= \det(\lambda \textbf{I}-\textbf{A}(G))$. An \emph{eigenvalue of } $G$ is any root of its characteristic polynomial. Since $\textbf{A}(G)$ is %real and
symmetric, all the eigenvalues of $G$ are real numbers. The
\emph{spectrum}
of $G$
is   the set of its eigenvalues together with their multiplicities.
%If $\textbf{A}_1$ and $\textbf{A}_2$ are adjacency matrices which arise from two different labeling of the graph $G$, then, for some permutation matrix $\textbf{P}$, $\textbf{A}_1= \textbf{P}^{-1}\textbf{A}_2\textbf{P}$.

%Two graphs which contain the same number of vertices connected in the same way are said to be \textit{isomorphic}.
 A \textit{walk} of length $k$ from vertex $x_0$ to $x_k $, denoted by $x_0\backsim x_1\backsim x_2\backsim \hdots\backsim x_{k-1}\backsim x_k$, is a sequence of vertices $x_0, x_1, x_2, \hdots, x_{k-1}, x_k$ such that the subsets $\{x_0 ,x_1\}, \{x_1 , x_2\}, \hdots, \{x_{k-1} , x_k\}$  are all edges. A \emph{path} is a walk with all vertices (and hence all edges) distinct. A \emph{cycle} is a path from a vertex back to itself (so the first and last vertices coincide). The cycle on $n$ vertices is denoted $C_n$.

\subsection {Circulant matrices}
We may recall   known facts about circulant and block circulant matrices. For more details, see for example \cite{Daviscir}, \cite{kra} and \cite{Tee}.

%\begin{definition}
A square matrix in which each row (after the first) has the elements of the previous
row shifted cyclically one place right is called a \emph{circulant} matrix. So, a circulant matrix has the form
%Philip R. Davis
%\cite{Daviscir}
%We denotes it as
$$
\mathbf{M}= circ(m_{0}, m_{1},\ldots , m_{n-1}) =
\left[
    \begin{array}{cccc}
      m_{0} & m_{1} & \ldots & m_{n-1} \\
      m_{n-1} & m_{0} & \ldots & m_{n-2} \\
      \vdots & \vdots & \ddots & \vdots \\
      m_{1} & m_{2} & \ldots & m_{0} \\
    \end{array}
  \right],
$$
%\end{definition}
that is, % if  $\mathbf {M}$, it holds that \  \  \
\  $
\mathbf{M}=[m_{i,j}]=[m_{0,j-i}],
$
where the subscripts are reduced  modulo $n$.

It is easy to see that $\mathbf{M}=circ(0,1,0,...,0,1)$ is the adjacency matrix of  $C_n$.   The following result can be found in \cite{kra}.
For concreteness, we express the
$n$th roots of unity as powers of $\xi= e^{\frac{2\pi\imath}{n}}$, where $\imath^2=-1 $.

%Let
%$\xi:=2\pi/n$.
%Then the $n$th roots of 1 are: $
%\xi_j = e^{\imath2\pij/n} = eij# = cos j# + i sin j# (j = 0, 1, 2, . . . , n − 1).$
%The following facts about   a
%circulant matrix $\textbf{M}=circ(m_{0}, m_{1},\ldots , m_{n-1})$ are well known (see, for example, \cite{Tee}).

\begin{proposition}\label{autovalores_circulante}
The eigenvalues of $\textbf{M}=circ(m_{0}, m_{1},\ldots , m_{n-1})$  are
%$=\left(
%    \begin{array}{ccccc}
%      m_{0} & m_{1} & \ldots & m_{n-2} & m_{n-1}\\
%      m_{n-1} & m_{0} & \ldots & m_{n-3} & m_{n-2}\\
%      \vdots & \vdots & \ddots & \vdots  & \vdots\\
%      m_{2} & m_{3} & \ldots & m_{0} & m_{1} \\
%      m_{1} & m_{2} & \ldots & m_{n-1} & m_{0} \\
%    \end{array}
%  \right)
%$
$$\lambda_{l}=m_0+\xi^{l}m_1 +...+ \xi^{(n-1)l} m_{n-1} , %\  \mbox{for  \ $l$, \  $0 \leq l \leq n-1$,}
$$
for \   $0 \leq l \leq n-1$,   each one with corresponding eigenvector
$
{\mathbf{v}}_{l}=[1, \xi^{l}, \xi^{2l}, ..., $
$\xi^{(n-1)l}]^{\top}.  %0\leq l\leq n-1.
$
% where  $\xi$ denotes a  primitive $n$-root of the unity.
 %The set $\{\emph{\textbf{v}}_{l}=[1, \xi^{l}, \xi^{2l}, ..., \xi^{(n-1)l}]^{T};\hbox{ }0\leq l\leq n-1\}$,
 % where $\xi$ denotes a primitive root of unity, is clearly linearly independent.
If
%$\xi$ denotes a primitive root of unity and
$\mathbf{V}=\left[   \begin{array}{ccccccc}
                   \mathbf{v}_{0}&|&\mathbf{v}_{1}&|&\ldots &|&\mathbf{v}_{n-1}\\
\end{array}
\right] $
is the matrix whose columns are the vectors $\mathbf{v}_l$ then $\mathbf{V'}= \frac{1}{\sqrt{n}}\mathbf{V}$ is a unitary matrix and $(\mathbf{{V'}})^{-1}\mathbf{M}\mathbf{V'}$ is a diagonal matrix, that is, $\mathbf{V'}$ (unitarily) diagonalizes $\mathbf{M}$.
For symmetric circulant matrices, the eigenvalues are real numbers.

%
%Furthermore, the (orthogonal)  matrix
%$\mathbf{V'}= \frac{1}{\sqrt{n}}\left(\begin{array}{cccc}
%                                                                                        \mathbf{v}_{0} & \mathbf{v}_{1} & \ldots & \mathbf{v}_{n-1} \\
%                                                                                      \end{array}
%                                                                                    \right)$,
%                                                                                  orthogonally  diagonalizes $\textbf{M}$. In fact, all  circulant $n \times n$ are (ortogonally) diagonalized by the same matrix $\mathbf{V'}$.
%
%                                                                                     %onde $\emph{\textbf{v}}_{l}=(1, \xi^{l}, \xi^{2l}, ..., \xi^{(n-1)l})^{t}\hbox{,  }0\leq l\leq n-1$ e $\xi$ indica a raiz n-esima da unidade, diagonaliza $\mathbf{M}$. De fato, $\mathbf{(V')^{-1}MV'} = \frac{1}{n}\mathbf{V^{-1}MV}$ que e uma matriz diagonal pelo Corol\'ario \ref{matrizes simul diago}.
\end{proposition}

%For concreteness, we express the
%$n$th roots of unity as powers of $\xi= e^{\frac{2\pi\imath}{n}}$, where $\imath^2=-1 $.

%We fix, once and for all, a choice of an $n$-primitive root of unity $\xi=e^{\frac{2\pi\imath}{n}}$ (ver Kra e Kalman para a melhor forma!!!!
%

%\begin{definition}
%Uma matriz $\mathbf{V}$ de ordem $n$ e dita unit\'aria quando $\mathbf{V V^*} = \mathbf{V^* V}=\mathbf{I}$, onde $\mathbf{I}$ e a matriz identidade de ordem $n$ e $\mathbf{V^*}$ e a  matriz adjunta de $\mathbf{V}$, ou seja, a matriz em que cada entrada $ij$ e o complexo conjugado da entrada $ji$ de $\mathbf{V}$, ou seja, $\mathbf{V^*}=\overline{\mathbf{V}^t}=\overline{\mathbf{V}}^t$.
%\end{definition}
%
%\begin{definition}
%Diz-se que duas matrizes $\mathbf{M}$ e $\mathbf{N}$ são simultaneamente diagonaliz\'aveis quando existe uma matriz invertível $\mathbf{V}$ tal que $\mathbf{V^{-1}MV}$ e $\mathbf{V^{-1}NV}$ são ambas matrizes diagonais.
%\end{definition}
%
%\begin{proposition}
%Todas as matrizes circulantes de ordem $n$ são simultaneamente diagonalizadas pela mesma matriz unit\'aria.
%\end{proposition}
%

\subsection {$I$-graphs}

%The I-graph I(n, j, k) is a graph with vertex set
%V (I(n, j, k)) = {u0, u1, . . . , un−1, v0, v1, . . . , vn−1}
%and edge set
%E(I(n, j, k)) = {uiui+j , uivi, vivi+k : i = 0, . . . , n − 1}.
%Since I(n, j, k) = I(n, k, j) we will usually assume that j ≤ k. Clearly G(n, k) = I(n, 1, k).
%Following the usual representation of these graphs where we draw vertices ui on one circle
%and vertices vi on another circle (with smaller radius),

Let fix $n,j,k \in \mathbb N$ with  $n\geq 3$, $1 \leq j, k<\frac{n}{2}$ and $j \leq k$. The $I$-graph $I(n,j,k)$ is the graph with vertex set $V(I(n,j,k))=\{a_i, b_i; 0\leq i \leq n-1\}$ and  edge set $E(I(n,j,k))=\{\{a_i,a_{i+j}\}, \{a_i,b_i\}, \{b_i, b_{i+k}\}; 0\leq i \leq n-1\}$, where addition is performed modulo $n$.

We may assume that $j \leq k$ since  $I(n, j, k) = I(n, k, j)$.  We consider $j,k <\frac{n}{2}$ because  $I (n, j, k)$,
$I (n, n -j, k)$ and $I (n, j, n-k)$
 are isomorphic  (neither $j=n/2$ nor $k=n/2$ furnish simple graphs). Thus,  $I(n,j,k)$ are cubic graphs (that is, regular graphs of degree 3) on $2n$ vertices.  The \emph{Petersen} \emph{graph} is $I(5,1,2)$. The class of $I$-graphs contains the well known class of
$G(n, k)= I(n,1,k)$, the so called (\cite{origemnome}) \emph{generalized} \emph{Petersen} \emph{graphs}, introduced in  \cite{origempetersen}. For more  on graphs $G(n, k)$ see  also \cite{alspach},\cite{petersen}, \cite{Fru}, \cite{Horvunit}, \cite{hiper} and references therein.

%We call an I-graph that is connected and not isomorphic to a generalized Petersen
%graph a proper I-graph.

%\begin{remark}
%\begin{itemize}
%  \item O grafo $G(n,k)$ possui ordem $2n$.
%  \item $G(n,k)$ e sempre um grafo cúbico, ou seja, todos os seus vertices têm grau 3.
%  \item Todo grafo de Petersen generalizado e conexo. Basta percebermos que para cada $i=0, \hdots, n-1$, os vertices $a_i$ formam um ciclo de tamanho $n$ e que cada vertice $b_i$ e adjacente a $a_i$.
%\end{itemize}
%\end{remark}

 %The paper \cite{petersen} determines the spectrum of generalized Petersen graphs.

 %Neste mesmo artigo, o autor afirma tambem que se $n$ e par e $k$ e ímpar então o grafo de Petersen generalizado $G(n,k)$ e bipartido.

%\begin{remark}
%Citaremos alguma características dos $I$-grafos.
%\begin{itemize}
%  \item Como $I(n,j,k) = I(n,k,j)$ vamos assumir que $j\leq k$
%  \item $I$-grafos são grafos 3-regulares.
%  \item e f\'acil perceber que $G(n,k) = I(n,1,k)$. Assim, todo grafo de Petersen generalizado e um $I$-grafo e portanto, todos os resultados a seguir tambem valem para os GPG. O grafo a figura \ref{fig:i1234} e um dos menores $I$-grafos que não são grafos de Petersen generalizados.
%\end{itemize}
%\end{remark}

 %Let $n,j,k \in \mathbb{N}$ be fixed, where   $n\geq 3$, $j<\frac{n}{2}$,  $k<\frac{n}{2}$ and $j \leq k$.
 %We denote  by $A(n,j)$ the  subgraph of $I(n,j,k)$  with
 The vertices $\{a_i; 0\leq i \leq n-1\}$ and edges $\{\{a_i,a_{i+j}\}; 0\leq i \leq n-1\}$ of the $I$-graph $I(n,j,k)$ form a subgraph which we  denote by $A(n,j)$.  On the other hand, the subgraph of $I(n,j,k)$ with vertices $\{b_i; 0\leq i \leq n-1\}$ and edges $\{\{b_i,b_{i+k}\}; 0\leq i \leq n-1\}$  will be denoted $B(n,k)$. The following proposition gives us an idea of the  structures of
 these two subgraphs respecting connectivity.

%\begin{figure}[H]
%\centering
%\includegraphics[width=3.8in,height=3.0in]{I_12_3_41.eps}
%\caption{$I(12,3,4)$}
%\label{fig:i1234}
%\end{figure}
%\end{example}

%\begin{example}\label{A,B}
%A seguir temos $A(12,3)$ e $B(12,4)$ que são subgrafos do grafo da Figura ~\ref{fig:i1234}
%\begin{figure}[H]
%\centering
%\includegraphics[width=5.5in,height=3.7in]{AB(12,3,4).eps} %\quad
%%\includegraphics[width=3.5in,height=2.4in]{B(12,4).eps}
%\caption{$A(12,3)$ e $B(12,4)$, subgrafos de $I(12,3,4)$.}
%\end{figure}
%\end{example}

%%%%%%%%%%%%%%%%%%%%%%%%%%%%%%%%%%%%%%
%%%%%%%%%%%%%% Componentes subgrafos %%%%%%%%%%%%

\begin{proposition}\label{subgrafodesconexo} The subgraph $B(n,k)$ has $d_{B}=\gcd(n,k)$ connected components, each one isomorphic to a cycle of length $\frac{n}{d_{B}}$.
Analogously, the subgraph $A(n,j)$ has $d_{A}=\gcd(n,j)$ connected components, each one isomorphic to a cycle of length  $\frac{n}{d_{A}}$.
\end{proposition}

\prova
Let $q=\frac{n}{d_{B}} \in {\mathbb{N}}$.
 For  $k=1$,  $\gcd(n,k)=1$ and then, by construction, the graph $B(n,k)$ has a connected component, namely, the cycle  $C_n$.
For  $k>1$, the vertex $b_0$ is adjacent to  $b_{0+k}$, which in turn is adjacent to  $b_{2k}$, and so on, until    vertex $b_{(q-1)k}$, which is adjacent to $b_{qk}$. But  $b_{qk}$ coincides with $b_{0}$, since
$qk\equiv0$(mod $n$): indeed, $qk=%k\frac{n}{d_B}=
n\frac{k}{d_B}$ with $\frac{k}{d_B}\in \mathbb{N}$. Then   two cases are possible, depending on $\gcd(n,k)=1$ or $\gcd(n,k)>1$.
     If $\gcd(n,k)=1$,  then $qk=%\frac{n}{1}\cdot k=
     n\cdot k$ and   $q=n$; since we go through all the vertices in  $B(n,k)$, we have a connected component which is a cycle of lenght $q$.
  In case $d_B = \gcd(n,k)>1$, it follows that $q=\frac{n}{d_B}<n$ and then,  the closed path on $q$ vertices from $b_0$ to $b_{qk}=b_0$   is covered.  Thus we obtain a cycle  isomorphic to $C_q$. Starting again from the first vertex to the right of $b_0$  not in the cycle  $b_0 ... b_{qk}$, say $b_l$, we obtain $b_l\backsim b_{l+k}\backsim b_{l+2k}\backsim ... \backsim b_{l+qk}=b_l$, a cycle of lenght  $q$ (like before). Since we have $n=qd_B$ vertices, we may  proceed  until we obtain  $d_B$ cycles of lenght $q$.

  Analogously, it can be shown that the subgraph $A(n,j)$ has  $d_{A}=\gcd(n,j)$ connected components, each one isomorphic to a cycle of length $\frac{n}{d_{A}}$.
\begin{flushright}
    $\square$
\end{flushright}
\begin{figure}[h]
\centering
\includegraphics[width=3.6in,height=2.8in]{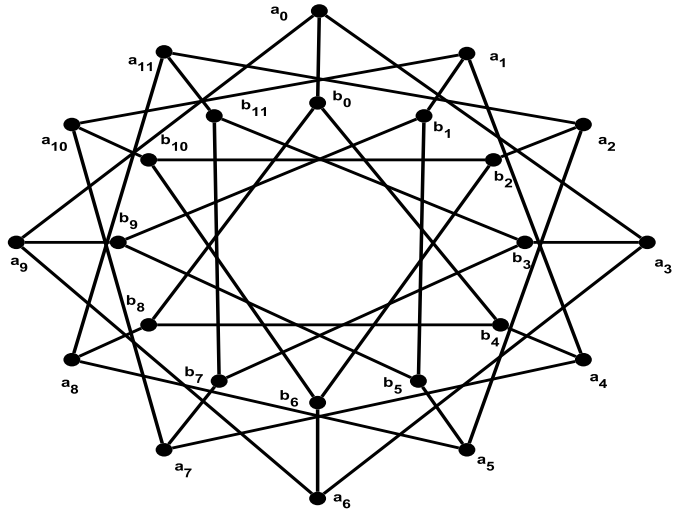}
\caption{$I(12,3,4)$}
\label{fig:i1234}
\end{figure}

 \begin{example}
 Figure~\ref{fig:i1234} depicts the $I$-graph $I(12,3,4)$  and the  vertices  $a_i$ and $b_i$ of  subgraphs $A(12,3)$ and $B(12,4)$, respectively.
 \end{example}

\vskip0.3cm

It follows from the above proposition  (or, alternatively,  from the $I$-graph  construction) that the adjacency matrices of the subgraphs $A(n,j)$ and  $B(n,k)$ are  circulant matrices.
% by the construction of  $A(n,j)$, its adjacency matrix $\textbf{A}(A(n,j))=[a_{rs}]$ satisfy
%$$
%a_{r,s}=\left\{
%         \begin{array}{ll}
%           1, & \hbox{if $s=r+j$ or $s=r-j$, reduced modulo n} \\
%           0, & \hbox{otherwise.}
%         \end{array}
%             \right.
%             $$
%% a_{0,s-r}=\left\{
%%  \begin{array}{ll}
%%    1, & \hbox{se $s-r=j$ ou $s-r=-j$ reduzidos módulo n } \\
%%    0, & \hbox{caso contr\'ario.}
%%  \end{array}
%%\right.
%$$=\left\{
%  \begin{array}{ll}
%    1, & \hbox{if $s-r=j$ or $s-r=-j$ reduced modulo $n$ } \\
%    0, & \hbox{otherwise,}
%  \end{array}
%\right.
%$$
%that is, $a_{r,s}=a_{0,s-r}$, for all $ r,s$,  $0\leq r,s\leq n-1$.    Furthermore,
Indeed,  $$\textbf{A}(A(n,j))=circ(\overbrace{0,...,0}^{j \ \hbox{{\small entries}}},1,0,...,0,1,\overbrace{0,...,0}^{j-1 \ \hbox{{\small {entries}}}})  \   \  \  \mbox{and}$$
%Similarly, it can be verified that the adjacency matrix of the subgraph $B(n,k)$ is
%$$
%b_{r,s}=\left\{
%         \begin{array}{ll}
%           1, & \hbox{se $s=r+k$ ou $s=r-k$, reduzidos módulo n} \\
%           0, & \hbox{caso contr\'ario.}
%         \end{array}
%       \right.
%$$
%Por outro lado, para $0\leq r,s\leq n-1$ vale
%$$
%b_{0,s-r}=\left\{
%  \begin{array}{ll}
%    1, & \hbox{se $s-r=k$ ou $s-r=-k$ reduzidos módulo n } \\
%    0, & \hbox{caso contr\'ario.}
%  \end{array}
%\right.
%$$
%ou seja,
%$$
%b_{0,s-r}=\left\{
%  \begin{array}{ll}
%    1, & \hbox{se $s=r+k$ ou $s=r-k$ reduzidos módulo n } \\
%    0, & \hbox{caso contr\'ario.}
%  \end{array}
%\right..
%$$
%
%Mas então, $b_{r,s}=b_{0,s-r}, \forall r,s$ com $0\leq r,s\leq n-1$. E assim, por definição, $\textbf{A}(B(n,k))$ é uma matriz circulante.
 $$\textbf{A}(B(n,k))
 =circ(\overbrace{0,...,0}^{k \ \hbox { {\small entries}}},1,0,...,0,1,\overbrace{0,...,0}^{k-1 \  \hbox{{\small entries}}})\,.$$
%\begin{figure}[h]
%\centering
%\includegraphics[width=3.6in,height=2.8in]{I_12_3_41.eps}
%\caption{$I(12,3,4)$}
%\label{fig:i1234}
%\end{figure}

%\begin{example}
% Figure~\ref{fig:i1234} depicts the $I$-graph $I(12,3,4)$ in its usual representation,  where we draw vertices $a_i$ on one circle
%and vertices $b_i$ on another circle (with smaller radius):
%%\begin{figure}[H]
%%\centering
%%\includegraphics[width=3.8in,height=3.0in]{I_12_3_41.eps}
%%\caption{$I(12,3,4)$}
%%\label{fig:i1234}
%%\end{figure}
%\end{example}
For instance, for   $I(12,3,4)$ it holds that
$\textbf{A}(A(12,3))= circ(0,0,0,1,0,0,0,$
$0,0,1,0,0)$  \ and \
$ \textbf{A}(B(12,4)) =circ(0,0,0,0,1,0,0,0,1,0,0,0).$
\vskip0.3cm

%%%%%%%%%%%%%%%%%%
For brevity,  we  denote
$\textbf{A}^{nj} = \textbf{A}(A(n,j))$ and $\textbf{B}^{nk} = \textbf{A}(B(n,k))$. The eigenvalues and eigenvectors of  these matrices may be determined by using
 Proposition \ref{autovalores_circulante}.

%A partir de agora, denotaremos a matriz de adjacência de $A(n,j)$ e por $\textbf{B}^{nk}$ a matriz de adjacência de $B(n,k)$.
%Como vimos na Proposição \ref{autovalores_circulante}, podemos determinar os autovalores e respectivos autovetores dessas matrizes. A partir desse resultado  temos os seguintes Corol\'arios.

\begin{lemma}\label{auto A} The eigenvalues of
 $\emph{\textbf{A}}^{nj}$
 %$=circ(\overbrace{0,...,0}^{j \hbox { {\small entries}}},1,0,...,0,1,\overbrace{0,...,0}^{j-1 \hbox { {\small entries}}})$
     are $\alpha_l=2\cos(\frac{2\pi jl}{n})$, % \ \ $0\leq l\leq n-1$,
     with correspond\-ing eigenvectors  $\emph{\textbf{v}}_l=\left(1, \xi^l, \xi^{2l},...,\xi^{(n-1)l}\right)^{\top}$, \ \ $0\leq l\leq n-1$.
 % ,where $\xi$ is a primitive $n$-root of unity.
  Analogously,   $\emph{\textbf{B}}^{nk}$ has eigenvalues
  %$=circ(\overbrace{0,...,0}^{k \hbox { {\small entries}}},1,0,...,0,1,\overbrace{0,...,0}^{k-1 \hbox { {\small entries}}})$,
   $\beta_l=2\cos(\frac{2\pi kl}{n})$, %\ \ $0\leq l\leq n-1$
    with corresponding eigenvectors $\emph{\textbf{v}}_l=(1, \xi^l, \xi^{2l},...,\xi^{(n-1)l})^{\top}$, \ \ $0\leq l\leq n-1$.
\end{lemma}

\prova
By Proposition \ref{autovalores_circulante}, the eigenvalues of  $\emph{\textbf{A}}^{nj}$ are  $\alpha_l=\xi^{jl}+\xi^{(n-j)l}$, \ \ $0\leq l\leq n-1$ , where
\begin{eqnarray*}
\alpha_l&=&\xi^{jl} + \xi^{(n-j)l}\\
		&=&\cos\left(\frac{2\pi jl}{n}\right)+i\sin\left(\frac{2\pi jl}{n}\right)+\cos\left(\frac{2\pi lj(n-j)}{n}\right)+i\sin\left(\frac{2\pi jk(n-j)}{n}\right)\\
		&=&\cos\left(\frac{2\pi jl}{n}\right)+i\sin\left(\frac{2\pi jl}{n}\right)+\cos\left(\frac{2\pi jl}{n}\right)-i\sin\left(\frac{2\pi l}{n}\right)\\
 &=&  2\cos\left(\frac{2\pi jl}{n}\right).
\end{eqnarray*}
 The other assertion  follows similarly.
\begin{flushright}
    $\square$
\end{flushright}

\section{Main result}

Let consider the $I$-graph $I(n,j,k)$ fixed, where $n,j,k \in \mathbb{N}$ are as before agreed.
%with  $n\geq 3$, $1 \leq j, k<\frac{n}{2}$ and $j \leq k$.
We label   the vertices of  $I(n,j,k)$ as follows: the vertices  $b_{0}, \ldots, b_{n-1}$ of $B(n,k)$ are designated, respectively,  as  $0,1,2, \ldots, s, \ldots, n-1$,  and the vertices $a_{0}, \ldots, a_{n-1}$ in $A(n,k)$ are labelled, respectively, as $n,n+1,n+2, \ldots, n+s, \ldots, 2n+1$.
From the adjacency rules in subgraphs $B(n,k)$ and $A(n,j)$  of $I(n,j,k)$, it follows respectively  that  $i\backsim (i+k)$ mod $n$ and  %$0\leq i\leq n-1$ and
 $(n+i)\backsim (n+i+j)$ mod $n$,  for $0\leq i\leq n-1$. Besides this, the adjacency rules between vertices of $A(n,j)$ and  $B(n,k)$ are given by
 $i \backsim (n+i)$, $0\leq i\leq n-1$.
This way, the adjacency matrix of $I(n,j,k)$ can be described as the following  block circulant matrix
$$
\mathbf{A}(I(n,j,k))=
\left(
  \begin{array}{cc}
   \mathbf{B}^{nk} & \mathbf{I}_{n} \\
    \mathbf{I}_{n} & \mathbf{A}^{nj} \\
  \end{array}
\right),
$$
where $\mathbf{I}_n$ is the identity matrix of order $n$.
 %and  $\emph{\textbf{A}}^{nj}$ and  $\emph{\textbf{B}}^{nk}$ are the adjacency matrix of the
%subgraphs $A(n,j)$ and $B(n,k)$ respectively.

\begin{theorem}\label{autovalores injk}
%Sejam $n,j,k \in \mathbb{N}$ fixos com $n\geq 3$, $j<\frac{n}{2}$ e $k<\frac{n}{2}$.
The eigenvalues of $I(n,j,k)$ have the form
$$
\lambda_{l} = \cos\left(\frac{2\pi jl}{n}\right)+\cos\left(\frac{2\pi kl}{n}\right)\pm \sqrt{\left(\cos\left(\frac{2\pi jl}{n}\right)-\cos\left(\frac{2\pi kl}{n}\right)\right)^2+1},$$
 for $ \ 0\leq l\leq n-1$.
\end{theorem}

\prova
For a fixed $l$, $0 \leq l \leq n-1$ let
$$
a_l=\frac{\beta_{l} - \alpha_{l} \pm \sqrt{(\alpha_{l} - \beta_{l})^{2}+4}}{2} \ \hbox{and} \ \lambda_l=\frac{\beta_{l} + \alpha_{l} \pm \sqrt{(\alpha_{l} + \beta_{l})^{2}-4(\alpha_l^2 \beta_l^2-1)}}{2},
$$
where $\alpha_l$ and $\beta_l$ are as given in Lemma \ref{auto A}.
We will show that   $\lambda_l$ is an eigenvalue of $I(n,j,k)$ with associated eigenvector $\mathbf{w}_l=\left(
  \begin{array}{c}
    a_l\mathbf{v}_l \\
    \mathbf{v}_l \\
  \end{array}
\right)$, where $\mathbf{v}_{l}=(1, \xi^{l}, \xi^{2l}, ..., \xi^{(n-1)l})^{\top}$. In fact, from Lemma \ref{auto A},   it follows
that
$$
\left(
  \begin{array}{cc}
    \mathbf{B}^{nk} & \mathbf{I}_n \\
    \mathbf{I}_n & \mathbf{A}^{nj} \\
  \end{array}
\right)
\left(
         \begin{array}{c}
           a_j \mathbf{v}_l \\
           \mathbf{v}_l \\
         \end{array}
       \right)=
\left(
  \begin{array}{c}
    a_l \mathbf{B}^{nk} \mathbf{v}_l + \mathbf{v}_l \\
    a_l \mathbf{v}_l + \mathbf{A}^{nj} \mathbf{v}_l \\
  \end{array}
\right)=
\left(
  \begin{array}{c}
    a_l \beta_l \mathbf{v}_l + \mathbf{v}_l \\
    a_l \mathbf{v}_l + \alpha_l \mathbf{v}_l \\
  \end{array}
\right)=
$$
\vspace{0.5cm}
$$
=\left(
  \begin{array}{c}
    (a_l \beta_l+1)\mathbf{v}_l \\
    (a_l+\alpha_l)\mathbf{v}_l \\
  \end{array}
\right)=
\left(
  \begin{array}{c}
    \left(\frac{\beta_{l} - \alpha_{l} \pm \sqrt{(\alpha_{l} - \beta_{l})^{2}+4}}{2}\cdot \beta_l +1\right)\mathbf{v}_l \\
    \left(\frac{\beta_{l} - \alpha_{l} \pm \sqrt{(\alpha_{l} - \beta_{l})^{2}+4}}{2}+\alpha_l\right)\mathbf{v}_l \\
  \end{array}
\right)=
$$
\vspace{0.5cm}
$$
=\left(
  \begin{array}{c}
    \left(\frac{\beta_l^2-\beta_l\alpha_l\pm \beta_l\sqrt{(\alpha_l-\beta_l)^2+4}}{2}+1\right)\mathbf{v}_l \\
    \left(\frac{\beta_l+\alpha_l\pm\sqrt{(\alpha_l-\beta_l)^2+4}}{2}\right)\mathbf{v}_l \\
  \end{array}
\right)=
\left(
  \begin{array}{c}
    \lambda_l a_l \mathbf{v}_l \\
    \lambda_l \mathbf{v}_l \\
  \end{array}
\right)=
\lambda_l \left(
  \begin{array}{c}
    a_l \mathbf{v}_l \\
    \mathbf{v}_l \\
  \end{array}
\right).
$$
Then $\lambda_l$ is an eigenvalue of $I(n,j,k)$. Furthermore,
\begin{eqnarray*}
\lambda_{l} &=& \frac{\beta_{l} + \alpha_{l} \pm \sqrt{(\alpha_{l} + \beta_{l})^{2}-4(\alpha_l \beta_l-1)}}{2}
			 \ = \ \frac{\beta_{l} + \alpha_{l} \pm \sqrt{(\alpha_{l} - \beta_{l})^{2}+4}}{2}=\\
			&=&\cos\left(\frac{2\pi jl}{n}\right)+\cos\left(\frac{2\pi kl}{n}\right) \pm \sqrt{\left(\cos\left(\frac{2\pi jl}{n}\right)-\cos\left(\frac{2\pi kl}{n}\right)\right)^2+1}.
\end{eqnarray*}
The assertion is proved.
 \begin{flushright}
    $\square$
\end{flushright}

%%%%%%%%%%%%%%%%%%%%%%%%%%%%%%%%%%%%%%%%%%%%%%%%%%%%%%%%%%%%%%%%
In the following remark we show how to obtain $\lambda_l$ and $a_l$ for all $l$, $0\leq l \leq n-1$.

\begin{remark}
In order to find the numbers  $a_l$ and $\lambda_l$ and the vectors  $\mathbf{w}_l$ in Theorem \ref{autovalores injk} we proceeded as follows:
For a fixed $l$, $0\leq l \leq n-1$ consider  $\mathbf{w}_{l}=\left(
                        \begin{array}{c}
                          a_{l}\mathbf{v}_{l} \\
                          \mathbf{v}_{l} \\
                        \end{array}
                      \right),$
where $\mathbf{v}_{l}=(1, \xi^{l}, \xi^{2l}, ..., \xi^{(n-1)l})^{\top}$ and $a_{l}$ is a constant to be determined.
We are searching for the solutions of the matrix equation
$$
\mathbf{A}(I(n,j,k))\mathbf{w}_{l} = \lambda\mathbf{w}_{l},
$$
 which is
$$
\left(
  \begin{array}{cc}
    \mathbf{B}^{nk} & \mathbf{I}_{n} \\
    \mathbf{I}_{n} & \mathbf{A}^{nj} \\
  \end{array}
\right).
%%%%%%%%%%%%%%%%%
%%%%%%%%%%%%%%%%%
\left(
                        \begin{array}{c}
                          a_{l}\mathbf{v}_{l} \\
                          \mathbf{v}_{l} \\
                        \end{array}
                      \right)=
%%%%%%%%%%%%%%
%%%%%%%%%%%%%%
\lambda \left(
                        \begin{array}{c}
                          a_{l}\mathbf{v}_{l} \\
                          \mathbf{v}_{l} \\
                        \end{array}
                      \right).
$$
It suffices to solve the system
$$\left\{
  \begin{array}{ll}
    \begin{array}{c}
      a_l \mathbf{B}^{nk} \mathbf{v}_l + \mathbf{v}_l=\lambda a_l \mathbf{v}_{l} \\
       a_l \mathbf{v}_l + \mathbf{A}^{nj} \mathbf{v}_l = \lambda \mathbf{v}_{l}
    \end{array}
  \end{array}
\right..
$$
From Lemma \ref{auto A} we have:
$$
\left\{
  \begin{array}{ll}
    \begin{array}{c}
      a_l \beta_{l} \mathbf{v}_l+ \mathbf{v}_l=\lambda a_l \mathbf{v}_{l} \\
       a_l \mathbf{v}_l + \alpha_{l} \mathbf{v}_l = \lambda \mathbf{v}_{l}
    \end{array}
  \end{array}
\right.
%%%%%%%%%%%%%%%%%%%%%%%
\Leftrightarrow
%%%%%%%%%%%%%%%%%%%%%%%
%\left\{
%  \begin{array}{ll}
%    \begin{array}{c}
%       a_l (\beta_{l} - \lambda) \mathbf{v}_l = - \mathbf{v}_{l} \\
%       (\alpha_{l} - \lambda) \mathbf{v}_l = - a_l \mathbf{v}_{l}
%    \end{array}
%  \end{array}
%\right.
%%%%%%%%%%%%%%%%%%%%%%%
%\Leftrightarrow
%%%%%%%%%%%%%%%%%%%%%%%
\left\{
  \begin{array}{ll}
    \begin{array}{c}
       a_l (\lambda - \beta_{l}) \mathbf{v}_l = \mathbf{v}_{l} \\
       (\lambda - \alpha_{l}) \mathbf{v}_l = a_l \mathbf{v}_{l}
    \end{array}
  \end{array}
\right.
%%%%%%%%%%%%%%%%%%%%%%%
\Leftrightarrow
%%%%%%%%%%%%%%%%%%%%%%%
\left\{
  \begin{array}{ll}
    \begin{array}{c}
       a_l (\lambda - \beta_{l}) = 1 \\
       (\lambda - \alpha_{l}) = a_l
    \end{array}
  \end{array}
\right.
$$
 By substituting the value $a_l$ into the first equation, we obtain
$$
\left\{
  \begin{array}{ll}
    \begin{array}{c}
        (\lambda - \beta_{l})(\lambda - \alpha_{l}) = 1 \qquad\mathbf{(*)}\\
       (\lambda - \alpha_{l}) = a_l\qquad\mathbf{(**)}
    \end{array}
  \end{array}
\right.
$$
Equation $\mathbf{(*)}$ furnishes
$$
\lambda_{l}^{+} = \frac{\beta_{l} + \alpha_{l} + \sqrt{(\alpha_{l} + \beta_{l})^{2}-4(\alpha_l \beta_l-1)}}{2}
 \  \  \  \   \emph{and} $$
 $$ \lambda_{l}^{-} = \frac{\beta_{l} + \alpha_{l} - \sqrt{(\alpha_{l} + \beta_{l})^{2}-4(\alpha_l \beta_l-1)}}{2}.
$$
Then, from equation (**) we have the following values for $a_l$:
$$
a_l^+=(\lambda_l^+ - \alpha_l) = \frac{\beta_{l} - \alpha_{l} + \sqrt{(\alpha_{l} - \beta_{l})^{2}+4}}{2}
 \   \  \  \mbox{and}  $$
$$a_l^-=(\lambda_l^- - \alpha_l) = \frac{\beta_{l} - \alpha_{l} - \sqrt{(\alpha_{l} - \beta_{l})^{2}+4}}{2}
$$
Thus,
\begin{eqnarray*}
\lambda_{l}^{+} &=& \frac{\beta_{l} + \alpha_{l} + \sqrt{(\alpha_{l} + \beta_{l})^{2}-4(\alpha_l^2 \beta_l^2-1)}}{2}
				 \ = \  \frac{\beta_{l} + \alpha_{l} + \sqrt{(\alpha_{l} - \beta_{l})^{2}+4}}{2} \ =\\
				&=&\cos\left(\frac{2\pi jl}{n}\right)+\cos\left(\frac{2\pi kl}{n}\right) + \sqrt{\left(\cos\left(\frac{2\pi jl}{n}\right)-\cos\left(\frac{2\pi kl}{n}\right)\right)^2+1},
\end{eqnarray*}
Then $\lambda_{l}^{+}$ is an eigenvalue of  $\mathbf{A}(I(n,j,k))$ with associated eigenvector  $\mathbf{w}^{+}_{l} = \left(
                        \begin{array}{c}
                          a_{l}^{+}\mathbf{v}_{l} \\
                          \mathbf{v}_{l} \\
                        \end{array}
                      \right).$
Similarly,
\begin{eqnarray*}
\lambda_{l}^{-} &=& \frac{\beta_{l} + \alpha_{l} - \sqrt{(\alpha_{l} + \beta_{l})^{2}-4(\alpha_l^2 \beta_l^2-1)}}{2}\\
				%&=&\frac{\beta_{l} + \alpha_{l} - \sqrt{(\alpha_{l} - \beta_{l})^{2}+4}}{2}\\
				&=&\cos\left(\frac{2\pi jl}{n}\right)+\cos\left(\frac{2\pi kl}{n}\right) - \sqrt{\left(\cos\left(\frac{2\pi jl}{n}\right)-\cos\left(\frac{2\pi kl}{n}\right)\right)^2+1},
\end{eqnarray*}
that is, $\lambda_{l}^{-}$ is an eigenvalue of  $\mathbf{A}(I(n,j,k))$ with corresponding eigenvector $\mathbf{w}^{-}_{l} = \left(
                        \begin{array}{c}
                          a_{l}^{-}\mathbf{v}_{l} \\
                          \mathbf{v}_{l} \\
                        \end{array}
                      \right).$
 Then, for all  $l$, $0\leq l \leq n-1$,  $\lambda^{+}_{l}, \lambda_{l}^{-}$ are the  $2n$ eigenvalues of  $I(n,j,k)$.
\end{remark}
\vskip 0.4cm

\begin{corollary}\label{autoraizes}
The  $2n$ eigenvalues of $I(n,j,k)$ are   the roots of the equations $(\lambda - \beta_{l})(\lambda - \alpha_{l}) = 1$, where $\alpha_l=2\cos(\frac{2\pi jl}{n})$ and $\beta_l=2\cos(\frac{2\pi kl}{n})$, for all  $l$, $0\leq l \leq n-1$.
\end{corollary}

\prova
It follows from equation (*).
\begin{flushright}
    $\square$
\end{flushright}

%%%%%%%%%%%%%%%%%%%%%%%%%%%%%%%%%%%%%%%%%%%%%%%%%%%%%%%
The spectrum  of a  $G(n,k)$ graph was determined in \cite{petersen}, where   the case $\gcd(n,k)=1$  was completely proved and the other case had its proof outlined.  Since  $G(n,k) =I(n,1,k)$,  we have:

\begin{corollary}[\cite{petersen}]\label{autopetersen}
The eigenvalues of  $G(n,k)$ are, for all $ 0\leq l\leq n-1$,
$$
\lambda_{l} = \cos\left(\frac{2\pi l}{n}\right)+\cos\left(\frac{2\pi kl}{n}\right)\pm \sqrt{\left(\cos\left(\frac{2\pi l}{n}\right)-\cos\left(\frac{2\pi kl}{n}\right)\right)^2+1}.$$
\end{corollary}

\section{Two structural properties of $I$-graphs through spectral approach}

%\begin{remark}
%No artigo \cite{petersen}, os autores utilizam tecnicas diferentes da apresentada aqui para demonstrar o Corol\'ario \ref{autopetersen}. A prova e feita em duas partes: a primeira considera grafos $G(n,k)$ tais que$\gcd(n,k)=1$ e a segunda considera os casos em que $\gcd(n,k)>1$. O primeiro caso e demonstrado de maneira similar ao feito na nossa demonstração do Teorema \ref{autovalores injk}. J\'a para o segundo caso, a demonstração parece estar equivocada, pois os autores afirmam que os autovetores serão da forma $\mathbf{w}_l=\left(a_1\mathbf{v'}_l, a_2\mathbf{v'}_l, \hdots, a_d\mathbf{v'}_l, \mathbf{v}_l\right)$, onde $\mathbf{v'}_l=\left( 1, \xi^l, \hdots\ \xi^{(n'-1)l}\right)^{\top}$, com $n'=\frac{n}{d}$. Esta demonstração sugere que os autores fazem uma nova rotulação na qual os ciclos $C_{\frac{n}{d}}$ do subgrafo $B(n,k)$ formam blocos de matrizes circulantes de tamanho $\frac{n}{d}$ o que mudaria tambem a matriz do subgrafo $A(n,1)$ e faria com que seus autovalores não fossem mais $\alpha_l$.
%
%O que não fica explicado na demonstração e como de alguma forma os autovetores $\mathbf{w}_l=\left(a_1\mathbf{v'}_l, a_2\mathbf{v'}_l, \hdots, a_d\mathbf{v'}_l, \mathbf{v}_l\right)$ coincidem com os que obtivemos em nossas contas, ou como os autovalores $\beta_l$ aparecem, j\'a que os autovalores devem ser uma combinação da forma apresentada no Corol\'ario \ref{autoraizes} (que tambem e proposto em \cite{petersen}).
%
%\end{remark}

The following are known structural properties of $I$-graphs, which have been demonstrated in \cite{igrafo} using  combinatorial and graph-theoretical ideas. After our Theorem \ref{autovalores injk}, we can prove them  by using an  spectral approach.

%Usaremos a Proposição \ref{proporegular} para provar o seguinte resultado, que nos diz quando um $I$-grafo e conexo.
%Este resultado foi provado em \cite{igrafo} usando tecnicas não espectrais.
To prove  the next theorem we use the  known fact (\cite{Cve1}, Theorem 3.23)  that a regular graph of degree $r$ is connected if and only if its (largest) eigenvalue $r$ has multiplicity equal to one.

\begin{theorem}[\cite{igrafo}]\label{igrafo conexo}
$I(n,j,k)$ is a connected graph if and only if $\gcd(n,j,k)=1$.
\end{theorem}
\prova
Consider the graph $I(n,j,k)$ and let $d= \gcd(n,j,k)$.
Suppose $d > 1$. For each $l=q\frac{n}{d} $, where \quad $q=0, ..., d-1$, it follows from Theorem  \ref{autovalores injk} that
\begin{eqnarray*}
\lambda_l    %&=&\cos\left(\frac{2\pi jqnd}{nd}\right) + \cos\left(\frac{2\pi kqnd}{nd}\right) + \sqrt{\left(\cos\left(\frac{2\pi jqnd}{nd}\right) - \cos\left(\frac{2\pi kqnd}{nd}\right)\right)^2+1}\\
                 &=&\cos\left(2\pi jq\right) + \cos\left(2\pi kq\right) + \sqrt{\left(\cos\left(2\pi jq\right) - \cos\left(2\pi kq\right)\right)^2+1}\\
			&=&\cos\left( 2q'\pi\right) + \cos\left(2q''\pi\right) + \sqrt{\left(\cos\left(2q'\pi\right) - \cos\left(2q''\pi\right)\right)^2+1}\\
			&=& 3,
\end{eqnarray*}
where $q'=j/d$ and $q''=k/d$.
Thus   3 is  an eigenvalue of multiplicity equal to $d$ and therefore,  $I(n,j,k)$ has $d$  connected components. Since $d > 1$, $I(n,j,k)$ is disconnected.

Conversely, suppose  $\gcd(n,j,k)=1$. Then $\gcd(n,j)=1$ or $\gcd(n,k)=1$ or $\gcd(j,k)=1$.
%In order for $I(n,j,k)$ to be connected, the multiplicity of the eigenvalue 3 must be one.
 We claim that, in any case, if $0\leq l \leq n-1$ and   $\lambda_l^+=3$ then $l=0$ (we note that $\lambda_l^ \leq 1$).
Indeed,  $\lambda_l^+=3$ if and only if
$$
3 = \cos\left(\frac{2\pi jl}{n}\right)+\cos\left(\frac{2\pi kl}{n}\right) + \sqrt{\left(\cos\left(\frac{2\pi jl}{n}\right)-\cos\left(\frac{2\pi kl}{n}\right)\right)^2+1}\quad \hbox{$\mathbf{(1)}$}
$$
Denoting $x=\cos\left(\frac{2\pi jl}{n}\right)$ and $y=\cos\left(\frac{2\pi kl}{n}\right)$,   we may solve the equation $3=x+y+\sqrt{(x-y)^2+1}$, that is, $2xy-3x-3y+4=0$, for $-1 \leq x, y \leq 1$.  We write
$
y=\frac{3x-4}{2x-3}=f(x), \quad x\neq \frac{3}{2}.
$
Since $f'(x)=-\frac{1}{(3-2x)^2} < 0$, the function $f$ decreases. Also $f(-1)=\frac{7}{5}$ and $f(1)=1$, and thus there is no $x \in (-1,1)$ such that $-1 \leq f(x) \leq 1$. Therefore, the only solution of the equation is $x=1$ and $y=1$.
Returning  to $\mathbf{(1)}$, the previous analysis leads us to
$$
\left\{
  \begin{array}{ll}
    (\cos(\frac{2\pi jl}{n})-\cos(\frac{2\pi kl}{n}))^2 = 0 \\
    \cos(\frac{2\pi jl}{n})=1 \\
    \cos(\frac{2\pi kl}{n})=1
  \end{array}
\right.
\hspace{0.2cm}\Leftrightarrow\hspace{0.3cm}
\left\{
  \begin{array}{ll}
    \frac{2\pi jl}{n}=2\pi r \\
    \frac{2\pi kl}{n}=2\pi r'
  \end{array}
\right.
\hspace{0.2cm}\Leftrightarrow\hspace{0.3cm}
\left\{
  \begin{array}{ll}
    jl=rn \\
    kl=r'n
  \end{array}
\right.
%\hspace{0.2cm}\Leftrightarrow\hspace{0.2cm}
%\left\{
%  \begin{array}{ll}
%    0\leq \frac{rn}{j} \leq n-1 \\
%    0\leq \frac{r'n}{k} \leq n-1
%  \end{array}
%\right.
$$
for $r,r'\in \mathbb{Z}$.
%So
%$
%\left\{
%  \begin{array}{ll}
%    0\leq \frac{rn}{j} \leq n-1 \\
%    0\leq \frac{r'n}{k} \leq n-1
%  \end{array}
%\right.
%$.
Thus $ 0\leq \frac{rn}{j} \leq n-1 $ and
$ 0\leq \frac{r'n}{k} \leq n-1$. We note that  $r$ is not a multiple of $j$, otherwise $\frac{rn}{j}\geq n$. Analogously, $r'$ is not a multiple of  $k$.

Let suppose $\gcd(n,j)=1$. From $n|jl$ %(we may recall that  $jl=rn$)
it follows that  $n|l$ and then,  $l<n$ implies  $l=0$.
Similarly, if  $\gcd(n,k)=1$ holds then  $l=0$.
For the case  $\gcd(j,k)=1$, let suppose   that $l\neq0$. Therefore $r$ and  and $r'$ are different from zero  and, from the  last system above, we obtain     $jr'=kr$.  Then the integer $j$ can be written as $j=\frac{kr}{r'}$ and, since $k$ is not a multiple of $j$, it follows that  $r'|k$, say, $k=tr'$ for an integer $t>1$. From  $jr'=tr'r$ we have $j=tr$, and thus $k$ and $j$ have a common divisor $t>1$, a contradiction. Then  $l=0$.
% For the case  $\gcd(j,k)=1$, let suppose by contradiction that $l\neq0$, which implies $r\neq 0$ and $r'\neq 0$. Thus, from the last system we obtain  $jr'=kr$. Then $j=\frac{kr}{r'}$. Notice that $r'\nmid r$. In fact, $\gcd(j,k)=1$ implies that $j$ is not a multiple of $k$. Since $j=\frac{kr}{r'} \in \mathbb{N}$, $r'|k$, i.e, $k=tr'$. Thus, $tr'l=r'n$ and form that $t=\frac{n}{l}$. So $jl=rtl$ and $j=tr$.
%Thus, $t|j$ and $t|k$ what is a contradiction.
%Then, $r=0$ or $r'=0$ hence $l=0$.
Therefore the multiplicity of the eigenvalue 3 is equal to one and the graph $I(n,j,k)$ is connected.
\begin{flushright}
    $\square$
\end{flushright}

In the sequence,  we use a known characterization of bipartite  graphs, which can be found in \cite{Cve1} (Theorem 3.11).

\begin{proposition} A graph containing at least one edge is bipartite if and only if its spectrum is
  symmetric with respect to zero.
\end{proposition}

\begin{theorem}[\cite{igrafo}]\label{igrafo bipartido}
The connected $I$-graph $I(n,j,k)$  is bipartite if and only if $n$ is even and $j$ and $k$ are odd.
\end{theorem}

\prova
Suppose that the connected graph $I(n,j,k)$ is  bipartite. As we have seen, $\lambda_0^+=3$ is an eigenvalue of $I(n,j,k)$ of multiplicity one. Since $I(n,j,k)$ is bipartite, exists $l$, $0\leq l \leq n-1$ such $\lambda_l^-=-3$
%Notice that does not exist $l$, $0\leq l \leq n-1$ such that $\lambda_l^+=-3$, since $  \lambda_l^+ \geq -1$, for $0\leq l \leq n-1$.)
(we may note that $  \lambda_l^+ \geq -1$, for all $0\leq l \leq n-1$).
For this fixed $l$ it holds that
%$$
%3+ \cos\left(\frac{2\pi jl}{n}\right) + \cos\left(\frac{2\pi kl}{n}\right) - \sqrt{\left(\cos\left(\frac{2\pi jl}{n}\right) - \cos\left(\frac{2\pi kl}{n}\right)\right)^2+1}=0,
%$$ hence,
%
$$
\cos\left(\frac{2\pi jl}{n}\right) + \cos\left(\frac{2\pi kl}{n}\right) - \sqrt{\left(\cos\left(\frac{2\pi jl}{n}\right) - \cos\left(\frac{2\pi kl}{n}\right)\right)^2+1}=-3.
$$
Reasoning as in the proof of Theorem \ref{igrafo conexo}, we conclude that the only possibility for the above equation in the
interval $[-1,1]$  is
$$
\left\{
  \begin{array}{ll}
    \cos\left(\frac{2\pi jl}{n}\right)=-1  \\
    \cos\left(\frac{2\pi kl}{n}\right)=-1 \\
    \left(\cos\left(\frac{2\pi jl}{n}\right) - \cos\left(\frac{2\pi kl}{n}\right)\right)^2=0
  \end{array}
\right.
%%%%%%%
 \  \  \  \Leftrightarrow  \  \
%%%%%%%%
\left\{
  \begin{array}{ll}
   \frac{2\pi j l}{n}=\pi + 2\pi r  \\
   \frac{2\pi k l}{n}=\pi + 2\pi r'
  \end{array}
\right.,
$$
for $r, r' \in \mathbb{Z}$. Thus
%$$
%\left\{
%  \begin{array}{ll}
%   \frac{2l}{n}=1+2r  \\
%   \frac{2k}{n}=1+2r'
%  \end{array}
%\right.
%$$
$\frac{2jl}{n}=1+2r$ and   $\frac{2kl}{n}=1+2r'$.
Therefore, $2jl=n(1+2r)$, and then $n$ is even. Also we have % $\frac{1+2r}{j}=\frac{1+2r'}{k}$, which implies
$j(1+2r)=k(1+2r')$, and so  $j$ and $k$ have  the same parity. But if both the integers are even, it follows that $d=\gcd(n,j,k)>1$,
%(since $n$ is also even),
 a contradiction according to Theorem \ref{igrafo conexo}, since the graph is connected. It follows that $n$ is even and $j$ and $k$ are odd.

For the converse,   suppose that $n$ is even and let $l\in\mathbb{N}$ be such that $0\leq l\leq \frac{n}{2}$. From trigonometry, for  odd integers $j$ and $k$   we have
$$
\cos\left(\frac{2\pi j\left(\frac{n}{2}-l\right)}{n}\right)= %\cos\left(\frac{2\pi jn}{2n}-\frac{2\pi jl}{n}\right)=
\cos(\pi j) \cos\left(\frac{2\pi jl}{n}\right)-\sin(\pi j) \sin\left(\frac{2\pi jl}{n}\right)= -\cos\left(\frac{2\pi jl}{n}\right) \ $$
%$$
%=\cos(\pi j) \cos\left(\frac{2\pi jl}{n}\right)-\sin(\pi j) \sin\left(\frac{2\pi jl}{n}\right)= -\cos\left(\frac{2\pi jl}{n}\right).
%$$
$$
\mbox{ and } \ \cos\left(\frac{2\pi k\left(\frac{n}{2}-l\right)}{n}\right)=-\cos\left(\frac{2\pi kl}{n}\right).
$$
Thus, from Theorem \ref{autovalores injk}, it follows that
%$$
%\lambda_l^+ + \lambda_{\frac{n}{2}-l}^-=\cos\left(\frac{2\pi jl}{n}\right)+\cos\left(\frac{2\pi kl}{n}\right)+\sqrt{\left(\cos\left(\frac{2\pi jl}{n}\right)-\cos\left(\frac{2\pi kl}{n}\right)\right)^2+1}\hbox{ +}
%$$
%$$
%+ \cos\left(\frac{2\pi j\left(\frac{n}{2}-l\right)}{n}\right)+\cos\left(\frac{2\pi k\left(\frac{n}{2}-l\right)}{n}\right)-\sqrt{\left(\cos\left(\frac{2\pi j\left(\frac{n}{2}-l\right)}{n}\right)-\cos\left(\frac{2\pi k\left(\frac{n}{2}-l\right)}{n}\right)\right)^2+1}=
%$$
$$
\lambda_l^+ + \lambda_{\frac{n}{2}-l}^-=\cos\left(\frac{2\pi jl}{n}\right)+\cos\left(\frac{2\pi kl}{n}\right)+\sqrt{\left(\cos\left(\frac{2\pi jl}{n}\right)-\cos\left(\frac{2\pi kl}{n}\right)\right)^2+1}
$$
$$
-\cos\left(\frac{2\pi jl}{n}\right)-\cos\left(\frac{2\pi kl}{n}\right)-\sqrt{\left(-\cos\left(\frac{2\pi jl}{n}\right)+\cos\left(\frac{2\pi kl}{n}\right)\right)^2+1}=0,
$$
%Fro trigonometry we obtain
%$$
%\cos\left(\frac{2\pi j\left(\frac{n}{2}-l\right)}{n}\right)= \cos\left(\frac{2\pi jn}{2n}-\frac{2\pi jl}{n}\right)=
%$$
%$$
%=\cos(\pi j) \cos\left(\frac{2\pi jl}{n}\right)-\sin(\pi j) \sin\left(\frac{2\pi jl}{n}\right)= -\cos\left(\frac{2\pi jl}{n}\right).
%$$
%Analogously, since $k$ is also odd
%$$
%\cos\left(\frac{2\pi k\left(\frac{n}{2}-l\right)}{n}\right)=-\cos\left(\frac{2\pi kl}{n}\right).
%$$
Also, if  $l$ satisfies $\frac{n}{2}+1\leq l \leq n-1$ then
$\displaystyle{
\lambda_l^+ + \lambda_{\frac{3n}{2}-l}^-=0,}$
%$$=\cos\left(\frac{2\pi jl}{n}\right)+\cos\left(\frac{2\pi kl}{n}\right)+\sqrt{\left(\cos\left(\frac{2\pi jl}{n}\right)-\cos\left(\frac{2\pi kl}{n}\right)\right)^2+1} \hbox{ +}
%$$
%$$
%+ \cos\left(\frac{2\pi j\left(\frac{3n}{2}-l\right)}{n}\right)+\cos\left(\frac{2\pi k\left(\frac{3n}{2}-l\right)}{n}\right)-\sqrt{\left(\cos\left(\frac{2\pi j\left(\frac{3n}{2}-l\right)}{n}\right)-\cos\left(\frac{2\pi k\left(\frac{3n}{2}-l\right)}{n}\right)\right)^2+1}=
%$$
%$$
%=\cos\left(\frac{2\pi jl}{n}\right)+\cos\left(\frac{2\pi kl}{n}\right)+\sqrt{\left(\cos\left(\frac{2\pi jl}{n}\right)-\cos\left(\frac{2\pi kl}{n}\right)\right)^2+1}
%$$
%$$
%-\cos\left(\frac{2\pi jl}{n}\right)-\cos\left(\frac{2\pi kl}{n}\right)+\sqrt{\left(-\cos\left(\frac{2\pi jl}{n}\right)+\cos\left(\frac{2\pi kl}{n}\right)\right)^2+1}=0,
%$$
since we have
$$
\cos\left(\frac{2\pi k\left(\frac{3n}{2}-l\right)}{n}\right)=-\cos\left(\frac{2\pi kl}{n}\right)$$
$$  \  \mbox{and} \  \
\cos\left(\frac{2\pi j\left(\frac{3n}{2}-l\right)}{n}\right)=-\cos\left(\frac{2\pi jl}{n}\right).
$$
%That implies
%$$
%\cos\left(\frac{2\pi j\left(\frac{3n}{2}-l\right)}{n}\right)
%=\cos(3\pi j) \cos\left(\frac{2\pi jl}{n}\right)-\sin(3\pi j) \sin\left(\frac{2\pi jl}{n}\right)= -\cos\left(\frac{2\pi jl}{n}\right).
%$$
%Hence
%$$
%\cos\left(\frac{2\pi k\left(\frac{3n}{2}-l\right)}{n}\right)=-\cos\left(\frac{2\pi kl}{n}\right).
%$$
%Therefore, the eigenvalues of $I(n,j,k)$ satisfy
%$$
%\left\{
%  \begin{array}{ll}
%    \lambda_l^+ + \lambda_{\frac{n}{2}-l}^-=0, \quad 0\leq l\leq \frac{n}{2}\\
%    \lambda_l^+ + \lambda_{\frac{3n}{2}-l}^-=0, \quad\frac{n}{2}+1\leq l \leq n-1
%  \end{array}
%\right.,
%$$
Thus, the spectrum is symmetric about  the origin and  $I(n,j,k)$ is a bipartite graph.
\begin{flushright}
    $\square$
\end{flushright}

We may note that the above result does not apply to disconnected graphs. For instance,
the disconnected graph $I(12,2,2)$ is bipartite, since it is formed from two copies of $I(6,1,1)$, which is bipartite.

%\begin{proposition}
%Se o grafo $I(n,j,k)$ e desconexo então ele consiste de $d=gcd(n,j,k)$ copias de $I\left(\frac{n}{d}, \frac{j}{d}, \frac{k}{d} \right)$.
%\end{proposition}

%O proximo resultado e tambem provado em \cite{igrafo} por meio de tecnicas não espectrais.
%\begin{theorem}[\cite{igrafo}]
%The connected $I$-graph $I(n,j,k)$  is bipartite if and only if $n$ is even and $j$ and $k$ are odd..
%\end{theorem}

The \emph{nullity} $\eta=\eta(G)$
 of a graph
$G$
is the multiplicity of
the number zero in its spectrum.  This graph-spectrum based invariant  has a noteworthy application in chemistry and also in mathematics as we can learn in the survey  \cite{surveynulidade} and the references therein.
 In \cite{rowlinson}, the author investigate  multiplicity of eigenvalues  of arbitrary cubic graphs. Applied to  an arbitrary $I$-graph, the main result of this paper concerning nullity asserts that

 \begin{proposition}[\cite{rowlinson}] The nullity $\eta$ of a   connected  $I$-graph $I(n,j,k)$ satisfies $\eta \leq n+1$.

 \end{proposition}

In the sequence, we determine the  nullity of $I$-graphs of the form $I(n,j,2j)$
by using  our Theorem \ref{autovalores injk}.

%\begin{lemma}
%Seja $I(n,j,2j)$ um $I$-grafo. Zero e autovalor de $I(n,j,2j)$ se e somente se existe $l$, $0\leq l \leq n-1$ tal que $3jl=n(r+1)$, para algum
%$r\not{\!\!\equiv} 2\emph{mod}$ 3 ou existe $l$, $0\leq l \leq n-1$ tal que $10jl=n(s+1)$, para algum $s\not{\!\!\equiv} 2$\emph{mod} 5.
%\end{lemma}

\begin{lemma}\label{lemanulidade} There exists $l$, $0\leq l \leq n-1$, such that $\lambda_l=0$  is an eigenvalue of the $I$-graph $I(n,j,2j)$ if and only if exists $l$, $0\leq l \leq n-1$, such that $3jl=n(1+3r)$  %,for $r \in\mathbb{Z}^+$,
 or $3jl=n(2+3r)$ %,  for $r \in\mathbb{Z}^+$,
 or
 $10jl=n(1+10r)$ %, $r \in\mathbb{Z}^+$,
  or
  $10jl=n(3+10r)$ %, $r \in\mathbb{Z}^+$,
   or
   $10jl=n(7+10r)$ %, $r \in\mathbb{Z}^+$,
   or
    $10jl=n(9+10r)$, where $r \in\mathbb{Z}^+$ in each case.
 \end{lemma}

\prova
From Theorem \ref{autovalores injk}, for $0\leq l \leq n-1$,  the eigenvalues of $I(n,j,2j)$ are of the form
\begin{eqnarray*}
\lambda_{l}^{\pm}&=& \cos\left(\frac{2\pi jl}{n}\right)+\cos\left(\frac{4\pi jl}{n}\right)\pm \sqrt{\left(\cos\left(\frac{2\pi jl}{n}\right)-\cos\left(\frac{4\pi jl}{n}\right)\right)^2+1}.
\end{eqnarray*}
 After  algebraic manipulations, we may give the above equations the form
$$
\lambda_{l}^{\pm}  \ = \ 2x^2+x-1 \pm \sqrt{\left(2x^2-x-1\right)^2+1},
$$
where $x=x(l)=\cos\left(\frac{2\pi jl}{n}\right)$ and $0\leq l \leq n-1$.
Thus, finding $l$, $0 \leq l \leq n-1$ such that $\lambda_l^{\pm}=0$ leads us to deal with the roots of the
equation $
8x^3-4x-1=0,
$
which are $
-\frac{1}{2}, \quad \frac{1+\sqrt{5}}{4}  \quad \hbox{and} \quad \frac{1-\sqrt{5}}{4}.
$
  It is straightforward to verify that  %$0\leq l\leq n-1$ such that
   $\lambda_l^{+}=0$ if and only %if exists such $l$ satisfying
$x(l) = -\frac{1}{2}$ \ or $x(l) = \frac{1-\sqrt{5}}{2}$ \ and   % there exists  $0\leq l\leq n-1$ such that
$\lambda_l^{-}=0$ if and only
%if exists such $l$ satisfying
$x(l) = \frac{1+\sqrt{5}}{2}$.

Let us analyse each case separately. Firstly,
$$
   \begin{array}{ll}
   \exists l \in \mathbb{Z}, 0\leq l\leq n-1, x(l)=\cos\left(\frac{2\pi jl}{n}\right)=-\frac{1}{2}\\
     \Leftrightarrow \ \exists l \in \mathbb{Z}, 0\leq l\leq n-1,  \ \frac{2\pi jl}{n}= \frac{2\pi}{3} +2\pi r, \ r \in \mathbb{Z}^+,   \vee
      \frac{2\pi jl}{n}= \frac{4\pi}{3} +2\pi r, \ r \in \mathbb{Z}^+,\\
       \Leftrightarrow \ \exists l \in \mathbb{Z}, 0\leq l\leq n-1,  \ 3jl= n(1+3r), \ r \in \mathbb{Z}^+,  \vee
     3jl = n(2+3r), \ r \in \mathbb{Z}^+\,.
     % \Leftrightarrow \ \exists l \in \mathbb{Z}, 0\leq l\leq n-1,  \ l= \frac{n(1+3r)}{3j}, \ r \in \mathbb{Z}^+  \vee
%      l = \frac{n(2+3r)}{3j}, \ r \in \mathbb{Z}^+\\
%     \Leftrightarrow \ \exists r \in \mathbb{Z}^+,  \ 3j|n(1+3r) \ \mbox{ and } 0 \leq r \leq  \frac{3jn-3j-n}{3n} \\
%     \vee  \  \
%     \ \exists r \in \mathbb{Z}^+,  \ 3j|n(2+3r) \ \mbox{ and } 0 \leq r \leq  \frac{3jn-3j-2n}{3n}.\\
   \end{array}
$$
Reasoning analogously, we obtain  the others two cases and complete the proof of our  assertion:
$$  \begin{array}{ll}
    \exists l \in \mathbb{Z}, \ 0\leq l\leq n-1,   \ \frac{2\pi jl}{n}= \frac{1-\sqrt{5}}{4}\\
   \Leftrightarrow \ \exists l \in \mathbb{Z},  \ 0\leq l\leq n-1,   \ \frac{2\pi jl}{n}= \frac{3\pi}{5} +2\pi r,  r \in \mathbb{Z}^+,    \vee    \frac{2\pi jl}{n}= \frac{7\pi}{5} +2\pi r, \ r \in {\mathbb{Z}}^+,\\
  \Leftrightarrow \  \exists l \in \mathbb{Z}, 0\leq l\leq n-1, \ 10jl=n(3+10r), r \in \mathbb{Z}^+,   \vee \
   10jl=n(7+10r), r \in {\mathbb{Z}}^+, \\
  % \Leftrightarrow \ \exists r \in \mathbb{Z}^+,  \ 10jl|n(3+10r)   \\
   %  \vee  \  \
     %\ \exists r \in \mathbb{Z}^+,  \ 10jl|n(7+10r)
%       \Leftrightarrow \ \exists r \in \mathbb{Z}^+,  \ 10j|n(3+10r) \ \mbox{ and } 0 \leq r \leq  \frac{10nj-10j-3n}{10n} \\
%     \vee  \  \
    %\ \exists r \in \mathbb{Z}^+,  \ 10j|n(7+10r) \ \mbox{ and } 0 \leq r \leq  \frac{10jn-10j-7n}{10n}  \   \   \
     \end{array}   $$
     and
$$\begin{array}{ll}
    \exists l \in \mathbb{Z},  \ 0\leq l\leq n-1,   \ \frac{2\pi jl}{n}= \frac{1+\sqrt{5}}{4}\\
   \Leftrightarrow \ \exists l \in \mathbb{Z}, 0\leq l\leq n-1,   \ \frac{2\pi jl}{n}= \frac{\pi}{5} +2\pi r, \ r \in \mathbb{Z}^+,   \vee    \frac{2\pi jl}{n}= \frac{9\pi}{5} +2\pi r, \ r \in \mathbb{Z}^+\\
      \Leftrightarrow  \exists l \in \mathbb{Z},   0\leq l\leq n-1,   10jl=n(1+10r),  r \in \mathbb{Z}^+, \vee   10jl=n(9+10r),  r  \in \mathbb{Z}^+.
      %\Leftrightarrow \ \exists r \in \mathbb{Z}^+,  \ 10jl|n(1+10r), \  \vee \  \exists r \in \mathbb{Z}^+,  \ 10jl|n(9+10r).
     %\ \exists r \in \mathbb{Z}^+,  \ 10j|n(9+10r) \ \mbox{ and } 0 \leq r \leq  \frac{10jn-10j-9n}{10n}  \   \   \
%       \Leftrightarrow \ \exists r \in \mathbb{Z}^+,  \ 10j|n(1+10r) \ \mbox{ and } 0 \leq r \leq  \frac{10nj-10j-n}{10n}%\\
%     \vee  \  \
%     \ \exists r \in \mathbb{Z}^+,  \ 10j|n(9+10r) \ \mbox{ and } 0 \leq r \leq  \frac{10jn-10j-9n}{10n}
   \end{array}
$$
\begin{flushright}
    $\square$
\end{flushright}

\begin{theorem}
The nullity  of the $I$-graph $I(n,j,2j)$ is  $\eta=|\mathcal{R}_1|+|\mathcal{R}_2|+|\mathcal{S}_1|+|\mathcal{S}_3|+|\mathcal{S}_7|+|\mathcal{S}_9|$, where %$\mathcal{R}_i$,$ \mathcal{S}_k$ are the sets
${\mathcal{R}}_i =\{r \in\ \mathbb{Z}  : 3j|n(i+3r) \ \wedge \ 0 \leq r \leq  j-1\}$
for $i=1,2$,
$\mathcal{S}_t =\{s \in\ \mathbb{Z}  : 10j|n(t+10s) \ \wedge \  0 \leq s \leq j-1  \}$, for $t=1,3,7,$ and
$\mathcal{S}_9 =\{s \in\ \mathbb{Z}  : 10j|n(9+10s) \ \wedge \  0 \leq s \leq j-1  \}$, if $n>10j$, or
$\mathcal{S}_9 =\{s \in\ \mathbb{Z}  : 10j|n(9+10s) \ \wedge \  0 \leq s \leq j-2  \}$, if $ n \leq 10j$.

\end{theorem}

\prova
 According Lemma \ref{lemanulidade},   to calculate the nullity we must counting the number of indices  $l$, $0\leq l \leq n-1$, satisfying any of its conditions. We will try to rewrite each of them in a simpler way.
   We begin noting that  from the above proof, provided $l$ is an integer, the equality $3jl=n(1+3r)$ means that $3j|n(1+r)$ ( then $l=\frac{n(1+r)}{3j}$) and also that  the integer $r$ satisfies  $0 \leq r \leq  \frac{3jn-3j-n}{3n} = j-\frac{j}{n} -\frac{1}{3}$, since $0 \leq l \leq n-1$. Furthermore, since  $ 2j<n/2$ as agreed from the beginning of the article, then $j-1 <j-(\frac{1}{4} +\frac{1}{3}) <j-(\frac{j}{n} +\frac{1}{3}) <j $; so, it suffices to have the integer $r$ satisfying  $0 \leq r \leq  j-1 = \lfloor  j-\frac{j}{n} -\frac{1}{3}\rfloor$.
%$$
%0\leq \frac{n(1+r)}{3j}\leq n-1 \quad\hbox{and then}\quad \quad $$
%$0\leq r \leq \frac{3jn-3j-n}{3n}$.
More precisely,  we just show that there exist $l, r \in \mathbb{Z}$ such that  $0 \leq l \leq n-1$   and $3jl=n(1+3r)$ (and in this case $l=\frac{n(1+3r)}{3j}$) if and only if there exists $r \in \mathbb{Z}$ such that  $3j|n(1+3r)$ and $0 \leq r \leq  j-1$.  This way, we obtain the elements of the set $\mathcal{R}_1$. By proceeding analogously, we may  determine the  sets
  $\mathcal{R}_2$, $\mathcal{S}_1$, $\mathcal{S}_3$ and  $\mathcal{S}_7$ and that $\mathcal{S}_9 =\{s \in \mathbb{Z} : 10j|n(9+10s) \wedge \ 0 \leq s \leq  j - \lfloor\j- \frac{j}{n}-\frac{9}{10}\rfloor\}$. Regarding the latter set, its description depends on having $n > 10j$ -- when it suffices to take $0 \leq s \leq  j-1$, -- or $n \leq 10j$, when we may request $0\leq s \leq  j-2$.
  We note that in case all of the sets are empty, it follows from  Lemma \ref{lemanulidade}  that  the $I$-graph $I(n,j,2j)$ has nullity equal to zero.
  \begin{flushright}
    $\square$
\end{flushright}

%That implies
%
%$$
%-1\leq r \leq (n-1)\frac{3j}{n}-1 \quad\hbox{so}\quad -1\leq r \leq \frac{3jn-3j-n}{n}
%$$
%provided $0 \leq l \leq n-1$.
%Analogously, from equation $10jl=n(2s+1)$ we must have $l=\frac{n(1+2s)}{10j}$, with  $10j|n(1+2s)$,  and
%$$
%0\leq \frac{n(1+2s)}{10j}\leq n-1 \quad\hbox{so}\quad 0\leq 1+2s \leq (n-1)\frac{10j}{n}
%$$
%Therefore,
%
%$$
%-1\leq 2s \leq (n-1)\frac{10j}{n}-1 \quad\hbox{so}\quad -1\leq s \leq \frac{10jn-10j-n}{2n}
%$$
%$-1\leq s \leq \frac{10jn-10j-n}{2n}$.
%From Lemma \ref{lemanulidade}, in case all   the sets $\cal{R}_i$ and $\cal{S}_k$ are empty, the I-graph $I(n,j,2j)$ has nullity equal to one.
%Therefore, the multiplicity of the eigenvalue zero is given by $\eta=|O_1|+|O_2|$.

\begin{example}
For instance,   $\eta(I(30,2,4))=4$, since  we have
$$ {\cal{R}}_1 = \left\{ r\in \mathbb{Z}:
 0 \leq r\leq 1 \wedge
  \frac{30(1+3r)}{6} \in \mathbb{Z} \right\}=\{0,1\}$$,
  $${\cal{R}}_2 = \left\{ r\in \mathbb{Z}; 0 \leq r\leq 1 \wedge \frac{30(2+3r)}{6} \in \mathbb{Z} \right\}=\{0,1\},$$
  $${\cal{S}}_1 = \left\{ r\in \mathbb{Z}; 0 \leq r\leq 1
  \wedge \frac{30(1+10r)}{20} \in \mathbb{Z} \right\}=\emptyset,$$
  $${\cal{S}}_3 = \left\{ r\in \mathbb{Z}; 0 \leq r\leq 1
  \wedge \frac{30(3+10r)}{20} \in \mathbb{Z} \right\}=\emptyset,$$
 $${\cal{S}}_7 = \left\{ r\in \mathbb{Z}; 0 \leq r\leq 1 \wedge \frac{30(7+10r)}{20} \in \mathbb{Z} \right\}=\emptyset$$   and
$${\cal{S}}_9 = \left\{ r\in \mathbb{Z}; 0 \leq r\leq 1 \wedge \frac{30(9+10r)}{20} \in \mathbb{Z} \right\}=\emptyset.$$ The indices are $5^+, 20^+, 10^+ $ and $25^+$ (see the proof of Lemma \ref{lemanulidade}).
 In turn,  $\eta(I(30,7,14))= 6$, since in this case,
$$ {\cal{R}}_1 = \left\{ r\in \mathbb{Z}:
 0 \leq r\leq 6 \wedge
  \frac{30(1+3r)}{21} \in \mathbb{Z} \right\}=\{2\},$$
 $${\cal{R}}_2 = \left\{ r\in \mathbb{Z}; 0 \leq r\leq 6 \wedge
   \frac{30(2+3r)}{21} \in \mathbb{Z} \right\}=\{4\},$$
 $${\cal{S}}_1 = \left\{ r\in \mathbb{Z};
  0 \leq r\leq 6   \wedge \frac{30(1+10r)}{70}
   \in \mathbb{Z} \right\}
  =\{4\},$$
 $${\cal{S}}_3 = \left\{ r\in \mathbb{Z}; 0 \leq r\leq 6
  \wedge \frac{30(3+10r)}{70} \in \mathbb{Z} \right\}=\{6\},$$
 $${\cal{S}}_7 = \left\{ r\in \mathbb{Z};    0 \leq r\leq 6
  \wedge  \frac{30(7+10r)}{70} \in \mathbb{Z} \right\}
 =\{0\}$$   and
 $${\cal{S}}_9 = \left\{ r\in \mathbb{Z}; 0 \leq r\leq 5 \wedge \frac{30(9+10r)}{70} \in \mathbb{Z} \right\}=\{4\}.$$ Here the indices are $10^+$, $20^+$, $9^+$, $ 27^+$, $3^-$ and $21^-$.

%\begin{eqnarray*}
%\cal{R}_1 &=& \left\{ r\in \mathbb{Z}; 0 \leq r\leq 1 \wedge \frac{30(1+3r}{6} \in \mathbb{Z} \right\}=\{0,1\},
%\cal{R}_1 &=& \left\{ r\in \mathbb{Z}; 0 \leq r\leq 1 \wedge \frac{30(1+3r}{6} \in \mathbb{Z} \right\}=\{0,1\},
% = \left\{ r\in \mathbb{Z}; 0\leq r\leq \frac{3}{2} \wedge 6|12(1+r)  \right\}= \{0,1\}\\
% &=&\left\{0,1,3,4\right\}
%\end{eqnarray*}
%e
%\begin{eqnarray*}
%O_2 &=& \left\{ s\in \mathbb{Z^+}; s\leq\frac{10jn-10j-n}{2n}, 10j|n(1+2s) \hbox{ e } s\not\equiv2\mod 5 \right\} \\
% &=& \left\{ s\in \mathbb{Z^+}; s\leq\frac{208}{24}, 20|12(1+2s) \hbox{ e } s\not\equiv2\mod 5 \right\}\\
% &=&\emptyset
%\end{eqnarray*}
\end{example}

\section{Conclusion}
In this paper, we completely determine the spectrum of an arbitrary $I$-graph and show how to use the theorem to obtain structural and spectral properties of these graphs.

%\section*{References}
%\begin{thebibliography}{1}
%\bibliography{Cybibfile}
%\bibliographystyle{plain}
%\bibliography{Cybfile}

\end{document}